\documentclass[12pt]{amsart}
\input{epsf}
\usepackage{amssymb,latexsym,cite}
\usepackage{graphicx}

\usepackage{epsfig}
\usepackage{subfigure}
\usepackage{caption}

\topskip  \newtheorem{theorem}{Theorem}

\usepackage{diagbox}
\usepackage{tabularx}
\usepackage{amsaddr}

\begin{document}

\pagenumbering{arabic}
\pagestyle{headings}

\title[ degenerate Fubini polynomials associated with random variables ]{Probabilistic degenerate Fubini polynomials associated with random variables}

\author{Rongrong  Xu$^{1}$,Taekyun  Kim$^{1,2,*}$, Dae San Kim$^{3,*}$,and Yuankui Ma$^{1,*}$}
\address{$^{1}$School of Science, Xi’an Technological University, Xi’an 710021, Shaanxi, China\\
$^{2}$Department of Mathematics, Kwangwoon University, Seoul 01897, Republic of Korea  \\
$^{3}$Department of Mathematics, Sogang University, Seoul 121-742, Republic of Korea }
\email{xurongrong0716@163.com; tkkim@kw.ac.kr; dskim@sogang.ac.kr; mayuankui@xatu.edu.cn }
\thanks{$^*$ Corresponding authors}

\begin{abstract}
Let $Y$ be a random variable such that the moment generating function of $Y$ exists in a neighborhood of the origin. The aim of this paper is to study probabilistic versions of the degenerate Fubini polynomials and the degenerate Fubini polynomials of order $r$, namely the probabilisitc degenerate Fubini polynomials associated with $Y$ and the probabilistic degenerate Fubini polynomials of order $r$ associated with $Y$. We derive some properties, explicit expressions, certain identities and recurrence relations for those polynomials. As special cases of $Y$, we treat the gamma random variable with parameters $\alpha,\beta > 0$, the Poisson random variable with parameter $\alpha >0$, and the Bernoulli random variable with probability of success $p$.
\end{abstract}

\maketitle

\section{Introduction}
In recent years, degenerate versions, $\lambda$-analogues and probabilistic versions of many special polynomials and numbers have been investigated by employing various methods such as  generating functions, combinatorial methods, umbral calculus, $p$-adic analysis, differential equations, probability, special functions, analytic number theory and operator theory (see [11-16, 18-21] and the references therein).\par
Let $Y$ be a random variable satisfying the moment condition (see \eqref{17}). The aim of this paper is to study probabilistic versions of the degenerate Fubini polynomials and the degenerate Fubini polynomials of order $r$, namely the probabilisitc degenerate Fubini polynomials associated with $Y$ and the probabilistic degenerate Fubini polynomials of order $r$ associated with $Y$. We derive some properties, explicit expressions, certain identities and recurrence relations for those polynomials and numbers. In addition, we consider the special cases that $Y$ is the gamma random variable with parameters $\alpha,\beta > 0$, the Poisson random variable with parameter $\alpha (>0)$, and the Bernoulli random variable with probability of success $p$. \par
The outline of this paper is as follows. In Section 1, we recall the degenerate exponentials, the degenerate Stirling numbers of the second kind ${n \brace k}_{\lambda}$, the degenerate Bell polynomials, the degenerate Fubini polynomials and the degenerate Fubini polynomials of order $r$. We remind the reader of Lah numbers and the partial Bell polynomials. Assume that $Y$ is a random variable such that the moment generating function of $Y$,\,\, $E[e^{tY}]=\sum_{n=0}^{\infty}\frac{t^{n}}{n!}E[Y^{n}], \quad (|t| <r)$, exists for some $r >0$. Let $(Y_{j})_{j\ge 1}$ be a sequence of mutually independent copies of the random variable $Y$, and let $S_{k}=Y_{1}+Y_{2}+\cdots+Y_{k},\,\, (k \ge 1)$,\,\, with \, $S_{0}=0$. Then we recall the probabilistic degenerate Stirling numbers of the second kind associated with $Y$ and the probabilistic degenerate Bell polynomials associated with $Y$, $\phi_{n,\lambda}^{Y}(x)$. Also, we remind the reader of the gamma random variable with parameters $\alpha,\beta > 0$. Section 2 is the main result of this paper. Let $(Y_{j})_{j \ge1},\,\, S_{k},\,\, (k=0,1,\dots)$ be as in the above. Then we first define the probabilistic degenerate Fubini polynomials associated with the random variable $Y$, $F_{n, \lambda}^{Y}(x)$. We derive for $F_{n, \lambda}^{Y}(x)$ an explicit expression in Theorem 2.1 and an expression as an infinite sum involving $E[(S_{k})_{n,\lambda}]$ in Theorem 2.2. In Theorem 2.3, when $Y \sim \Gamma(1,1)$, we find an expression for $F_{n, \lambda}^{Y}(x)$ in terms of Lah numbers and Stirling numbers of the first kind. We obtain a representation of $F_{n, \lambda}^{Y}(x)$ as an integral over $(0,\infty)$ of the integrand involving $\phi_{n,\lambda}^{Y}(x)$ and its generalization in Theorem 2.14. In Theorem 2.5, we express the probabilistic degenerate Fubini numbers associated with $Y$, $F_{n,\lambda}^{Y}=F_{n,\lambda}^{Y}(1)$, as a finite sum involving the partial Bell polynomials. Then we introduce the probabilistic degenerate Fubini polynomials of order $r$ associated with $Y$ and deduce an explicit expression for them in Theorem 2.6. We obtain a recurrence relation for $F_{n, \lambda}^{Y}(x)$ in Theorem 2.7, and another one in Theorem 2.8 together with its generalization in Theorem 2.15. In Theorem 2.9, the $r$th derivative of $F_{n, \lambda}^{Y}(x)$ is expressed in terms of $F_{i,\lambda}^{(r+1,Y)}(x)$. We get the identity $\frac{1}{1-x}F_{n,\lambda}^{Y}\bigg(\frac{x}{1-x}\bigg)=\sum_{i=0}^{\infty}E[(S_{i})_{n,\lambda}]x^{i}$ in Theorem 2.10 and its generalization in Theorem 2.13. In Theorem 2.11, when $Y$ is the Poisson random variable with parameter $\alpha$, we express $F_{n, \lambda}^{Y}(x)$ in terms of the Fubini polynomials $F_{i}(x)$ and ${n \brace i}_{\lambda}$. In Theorem 12, when $Y$ is the Poisson random variable with parameter $\alpha$, we show $\frac{1}{1-x}F_{n,\lambda}^{Y}\bigg(\frac{x}{1-x}\bigg)=\sum_{k=0}^{\infty}\phi_{n,\lambda}(k\alpha)x^{k}$. Finally, we show that $F_{n,\lambda}^{Y}(x)=F_{n,\lambda}(xp)$ if $Y$ is the Bernoulli random variable with probability of success $p$. For the rest of this section, we recall the facts that are needed throughout this paper.

\vspace{0.1in}

For any $\lambda\in\mathbb{R}$, the degenerate exponentials are defined by
\begin{equation}
e_{\lambda}^{x}(t)=\sum_{k=0}^{\infty}\frac{(x)_{k,\lambda}}{k!}t^{k},\quad e_{\lambda}(t)=e_{\lambda}^{1}(t),\quad (\mathrm{see}\ [6-21]),\label{1}
\end{equation}
where
\begin{equation}
(x)_{0,\lambda}=1,\quad (x)_{n,.\lambda}=x(x-\lambda)\cdots\big(x-(n-1)\lambda\big),\quad (n\ge 1).\label{2}
\end{equation}
Note that
\begin{displaymath}
	\lim_{\lambda\rightarrow 0}e_{\lambda}^{x}(t)=e^{xt}.
\end{displaymath}
The Stirling numbers of the first kind are defined by
\begin{equation}
(x)_{n}=\sum_{k=0}^{n}S_{1}(n,k)x^{k},\quad (n\ge 0),\quad (\mathrm{see}\ [1-3,5,24]), \label{3}
\end{equation}
where
\begin{displaymath}
	(x)_{0}=1,\quad (x)_{n}=x(x-1)\cdots(x-n+1),\quad (n\ge 1).
\end{displaymath}
The Lah numbers are defined by
\begin{equation}
\langle x\rangle_{n}=\sum_{k=0}^{n}L(n,k)(x)_{k},\quad (n\ge 0),\quad (\mathrm{see}\ [5-24]),\label{4}
\end{equation}
where
\begin{displaymath}
\langle x\rangle_{0}=1, \quad	\langle x\rangle_{n}=x(x+1)\cdots(x+n-1),\quad (n \ge 1).
\end{displaymath}
By \eqref{4}, we easily get
\begin{equation}
L(n,k)=\frac{n!}{k!}\binom{n-1}{k-1},\quad (n\ge k\ge 0),\quad (\mathrm{see}\ [5-24]). \label{5}
\end{equation}
In [13], the degenerate Stirling numbers of the second kind are defined by
\begin{equation}
	(x)_{n,\lambda}=\sum_{k=0}^{n}{n \brace k}_{\lambda}(x)_{k},\quad (n\ge 0). \label{6}
\end{equation}
It is well known that the degenerate Bell polynomials are defined by
\begin{equation}
e^{x(e_{\lambda}(t)-1)}=\sum_{n=0}^{\infty}\phi_{n,\lambda}(x)\frac{t^{n}}{n!},\quad (\mathrm{see}\ [12-14]).\label{7}	
\end{equation}
Thus, by \eqref{6} and \eqref{7}, we get
\begin{equation}
\phi_{n,\lambda}(x)=\sum_{k=0}^{n}{n \brace k}_{\lambda}x^{k},\quad (n\ge 0),\quad (\mathrm{see}\ [12,17,21]).\label{8}
\end{equation}
The degenerate Fubini polynomials are defined by
\begin{equation}
F_{n,\lambda}(x)=\sum_{k=0}^{n}{n\brace k}_{\lambda}k!	x^{k},\quad (n\ge 0),\quad (\mathrm{see}\ [17-19,26]).\label{9}
\end{equation}
Thus, by \eqref{9}, we get
\begin{equation}
\frac{1}{1-x(e_{\lambda}(t)-1)}=\sum_{n=0}^{\infty}F_{n,\lambda}(x)\frac{t^{n}}{n!},\quad (\mathrm{see} [4,19,21,26]). \label{10}	
\end{equation}
From \eqref{10}, we note that
\begin{equation}
\frac{1}{1-x}F_{n,\lambda}\bigg(\frac{x}{1-x}\bigg)=\bigg(x\frac{d}{dx}\bigg)_{n,\lambda}\frac{1}{1-x}=\sum_{k=0}^{\infty}(k)_{n,\lambda}x^{k},\quad (\mathrm{see}\ [18]). \label{11}
\end{equation}
For $r\in\mathbb{N}$, the degenerate Fubini polynomials of order $r$ are defined by
\begin{equation}
\bigg(\frac{1}{1-y(e_{\lambda}(t)-1)}\bigg)^{r}=\sum_{n=0}^{\infty}F_{n,\lambda}^{(r)}(y)\frac{t^{n}}{n!},\quad(\mathrm{see}\ [8,9,16]).\label{12}
\end{equation}
Thus, by \eqref{12}, we get
\begin{equation}
F_{n,\lambda}^{(r)}(y)=\sum_{k=0}^{n}\binom{k+r-1}{k}y^{k}{n\brace k}_{\lambda}k!,\quad (\mathrm{see}\ [8,18,19,22,23]). \label{13}	
\end{equation}
From \eqref{12}, we have
\begin{equation}
	\bigg(\frac{1}{1-x}\bigg)^{r+1}F_{n,\lambda}^{(r+1)}\bigg(\frac{x}{1-x}\bigg)=\bigg(x\frac{d}{dx}\bigg)_{n,\lambda}\bigg(\frac{1}{1-x}\bigg)^{r+1}=\sum_{k=0}^{\infty}\binom{k+r}{r}(k)_{n,\lambda}x^{k}, \label{14}
\end{equation}
where $n,r$ are nonnegative integers. \par
For any integer $k\ge 0$, the partial Bell polynomials are given by
\begin{equation}
	\frac{1}{k!}\bigg(\sum_{i=1}^{\infty}x_{i}\frac{t^{i}}{i!}\bigg)^{k}=\sum_{n=k}^{\infty}B_{n,k}(x_{1},x_{2},\dots,x_{n-k+1})\frac{t^{n}}{n!},\quad (\mathrm{see}\ [5]) \label{15}
\end{equation}
where
\begin{equation}
\begin{aligned}
&B_{n,k}(x_{1},x_{2},\dots,x_{n-k+1})\\
&=\sum_{\substack{l_{1}+l_{2}+\cdots+l_{n-k+1}=k\\ l_{1}+2l_{2}+\cdots+(n-k+1)l_{n-k+1}=n}}\frac{n!}{l_{1}l_{2}!\cdots l_{n-k+1}!}\bigg(\frac{x_{1}}{1!}\bigg)^{l_{1}} \bigg(\frac{x_{2}}{2!}\bigg)^{l_{2}}\cdots \bigg(\frac{x_{n-k+1}}{(n-k+1)!}\bigg)^{l_{n-k+1}}.
\end{aligned}\label{16}
\end{equation}
Let $Y$ be a random variable such that the moment generating function of $Y$
\begin{equation}
E[e^{tY}]=\sum_{n=0}^{\infty}E[Y^{n}]\frac{t^{n}}{n!},\quad (|t|<r)\quad \textrm{exists for some $r>0$}.\label{17}
\end{equation}
Assume that $(Y_{j})_{j\ge 1}$ is a sequence of mutually independent copies of $Y$ and $S_{k}=Y_{1}+Y_{2}+\cdots+Y_{k},\ (k\ge 1)$ with $S_{0}=0$.\par
The probabilistic degenerate Stirling numbers of the second kind associated with random variable $Y$ are defined by
\begin{equation}
{n \brace k}_{Y,\lambda}=\frac{1}{k!}\sum_{j=0}^{k}\binom{k}{j}(-1)^{k-j}E\big[(S_{j})_{n,\lambda}\big],\quad (n\ge k\ge 0),\quad (\mathrm{see}\ [15]). \label{18}
\end{equation}
By binomial inversion, the equation \eqref{18} is equivalent to
\begin{equation}
E\big[(S_{k})_{n,\lambda}\big]=\sum_{j=0}^{k}\binom{k}{j}j!{n \brace j}_{Y,\lambda},\quad (\mathrm{see}\ [15]).\label{19}
\end{equation}
From \eqref{18}, we note that
\begin{equation}
\frac{1}{k!}\Big(E[e_{\lambda}^{Y}(t)]-1\Big)^{k}=\sum_{n=k}^{\infty}{n \brace k}_{Y,\lambda}	\frac{t^{n}}{n!},\quad (k\ge 0),\quad (\mathrm{see}\ [15]). \label{20}
\end{equation}
In view of \eqref{8}, the probabilistic degenerate Bell polynomials associated with $Y$ are defined by
\begin{equation}
\phi_{n,\lambda}^{Y}(x)=\sum_{k=0}^{n}{n\brace k}_{Y,\lambda}x^{k},\quad (n\ge 0),\quad (\mathrm{see}\ [15,20]).\label{21}	
\end{equation}
When $Y=1$, we have $\phi_{n,\lambda}^{Y}(x)=\phi_{n,\lambda}(x)$.\par
By \eqref{21}, we get
\begin{equation}
e^{x(E[e_{\lambda}^{Y}(t)]-1}=\sum_{n=0}^{\infty}\phi_{n,\lambda}^{Y}(x)\frac{t^{n}}{n!},\quad (\mathrm{see}\ [15]).\label{22}	
\end{equation}
We recall that $Y$ is the gamma random variable with parameter $\alpha,\beta>0$ if probability density function of $Y$ is given by
\begin{equation*}
f(x)=\left\{\begin{array}{cc}
\frac{\beta}{\Gamma(\alpha)}e^{-\beta x}(\beta x)^{\alpha-1}, & \textrm{if $x\ge 0$},\\
0, & \textrm{if $x<0$,}
\end{array}\right.,\quad (\mathrm{see}\ [3,25,27]),
\end{equation*}
which is denoted by $Y\sim \Gamma(\alpha,\beta)$. \par

\section{Probabilistic degenerate Fubini polynomials associated with random variables}
Let $(Y_{k})_{k\ge 1}$ be a sequence of mutually independent copies of random variable $Y$, and let
\begin{displaymath}
	S_{0}=0,\quad S_{k}=Y_{1}+Y_{2}+\cdots+Y_{k},\quad (k\in\mathbb{N}).
\end{displaymath}
Now, we consider the {{\it probabilistic degenerate Fubini polynomials associated with random variable $Y$}} which are given by
\begin{equation}
\frac{1}{1-x\big(E[e_{\lambda}^{Y}(t)]-1\big)}=\sum_{n=0}^{\infty}F_{n,\lambda}^{Y}(x)\frac{t^{n}}{n!}.\label{23}
\end{equation}
For $Y=1$, we have $F_{n,\lambda}^{Y}(x)=F_{n,\lambda}(x),\ (n\ge 0)$. When $x=1$, $F_{n,\lambda}^{Y}=F_{n,\lambda}^{Y}(1)$ are called the {\it{probabilistic degenerate Fubini numbers associated with random variable $Y$}}. \par
From \eqref{23}, we note that
\begin{align}
&\frac{1}{1-x\big(E[e_{\lambda}^{Y}(t)]-1\big)}=\sum_{k=0}^{\infty}x^{k}\Big(E[e_{\lambda}^{Y}(t)]-1\Big)^{k} \label{24} \\	
&=\sum_{k=0}^{\infty}x^{k}k!\frac{1}{k!}\Big(E\big[e_{\lambda}^{Y}(t)\big]-1\Big)^{k}=\sum_{k=0}^{\infty}x^{k}k!\sum_{n=k}^{\infty}{n\brace k}_{Y,\lambda}\frac{t^{n}}{n!} \nonumber \\
&=\sum_{n=0}^{\infty}\sum_{k=0}^{n}x^{k}k!{n \brace k}_{Y,\lambda}\frac{t^{n}}{n!}. \nonumber
\end{align}
Therefore, by comparing the coefficients on both sides of \eqref{24}, we obtain the following theorem.
\begin{theorem}
For $n\ge 0$, we have
\begin{displaymath}
F_{n,\lambda}^{Y}(x)=\sum_{k=0}^{n}{n\brace k}_{Y,\lambda}k!x^{k}.
\end{displaymath}
\end{theorem}
By \eqref{24}, we get
\begin{align}
&\sum_{n=0}^{\infty}F_{n,\lambda}^{Y}(x)\frac{t^{n}}{n!}=\frac{1}{1-x(E[e_{\lambda}^{Y}(t)]-1)}=\frac{1}{1+x-xE[e_{\lambda}^{Y}(t)]} \label{25} \\
& =\frac{1}{1+x}\frac{1}{1-\frac{x}{1+x}E[e_{\lambda}^{Y}(t)]}=\frac{1}{1+x}\sum_{k=0}^{\infty}\bigg(\frac{x}{1+x}\bigg)^{k}\Big(E[e_{\lambda}^{Y}(t)]\Big)^{k} \nonumber \\
&=\sum_{n=0}^{\infty}\frac{1}{1+x}\sum_{k=0}^{\infty}\bigg(\frac{x}{1+x}\bigg)^{k}E\big[(Y_{1}+Y_{2}+\cdots+Y_{k})_{n,\lambda}\big]\frac{t^{n}}{n!}\nonumber \\
& =\sum_{n=0}^{\infty}\frac{1}{1+x}\sum_{k=0}^{\infty}\bigg(\frac{x}{1+x}\bigg)^{k}E[(S_{k})_{n,\lambda}]\frac{t^{n}}{n!}. \nonumber
\end{align}
Therefore, by \eqref{25}, we obtain the following theorem.
\begin{theorem}
For $n\ge 0$, we have
\begin{displaymath}
F_{n,\lambda}^{Y}(x)=\frac{1}{1+x}\sum_{k=0}^{\infty}\bigg(\frac{x}{1+x}\bigg)^{k}E[(S_{k})_{n,\lambda}].
\end{displaymath}
In particular, for $Y=1$, we have
\begin{displaymath}
F_{n,\lambda}^{Y}(x)=\frac{1}{1+x}\sum_{k=0}^{\infty}\bigg(\frac{x}{1+x}\bigg)^{k}(k)_{n,\lambda}.
\end{displaymath}
\end{theorem}
Let $Y\sim\Gamma(1,1)$. Then we have
\begin{align}
&\sum_{n=0}^{\infty}F_{n,\lambda}^{Y}(x)\frac{t^{n}}{n!}=\frac{1}{1-x(E[e_{\lambda}^{Y}(t)]-1)}=\sum_{k=0}^{\infty}x^{k}\Big(E[e_{\lambda}^{Y}(t)]-1\Big)^{k} \label{26} \\
&=\sum_{k=0}^{\infty}x^{k}\bigg(\int_{0}^{\infty}e_{\lambda}^{y}(t)e^{-y}dy-1\bigg)^{k}=\sum_{k=0}^{\infty}x^{k}\bigg(\int_{0}^{\infty}e^{y(\frac{1}{\lambda}\log (1+\lambda t)-1)}dy-1\bigg)^{k} \nonumber \\
&=\sum_{k=0}^{\infty}k!x^{k}\frac{1}{k!}\bigg(\frac{\frac{1}{\lambda}\log(1+\lambda t)}{1-\frac{1}{\lambda}\log(1+\lambda t)}\bigg)^{k}\nonumber \\
&=\sum_{k=0}^{\infty}k!x^{k}\sum_{l=k}^{\infty}L(l,k)\frac{1}{l!}\bigg(\frac{1}{\lambda}\log(1+\lambda t)\bigg)^{l}\nonumber \\
&=\sum_{l=0}^{\infty}\sum_{k=0}^{l}k!x^{k}L(l,k)\sum_{n=l}^{\infty}\lambda^{n-l}S_{1}(n,l)\frac{t^{n}}{n!}\nonumber \\
&=\sum_{n=0}^{\infty}\sum_{l=0}^{n}\sum_{k=0}^{l}k!x^{k}L(l,k)\lambda^{n-l}S_{1}(n,l)\frac{t^{n}}{n!}, \nonumber
\end{align}
where $S_{1}(n,l)$ are the Stirling numbers of the first kind. Here we should observe that, for all $t$ with $|t|$ small, we have
\begin{equation*}
\big| \frac{1}{\lambda} \log (1+ \lambda t) \big| <1,
\end{equation*}
since $|\frac{\log(1+x)}{x}|$ is bounded on $(0, \infty)$.
Therefore, by comparing the coefficients on both sides of \eqref{26}, we obtain the following theorem.
\begin{theorem}
Let $Y\sim\Gamma(1,1)$. Then we have
\begin{displaymath}
F_{n,\lambda}^{Y}(x)=\sum_{l=0}^{n}\sum_{k=0}^{l}k!\lambda^{n-l}L(l,k)S_{1}(n,l)x^{k},\quad (n\ge 0).
\end{displaymath}
\end{theorem}
Now, we observe from \eqref{21} and Theorem 2.1 that
\begin{align}
\int_{0}^{\infty}\phi_{n,\lambda}^{Y}(xy)e^{-y}dy &=\sum_{k=0}^{n}{n\brace k}_{Y,\lambda}x^{k}\int_{0}^{\infty}y^{k}e^{-y}dy \label{27} \\
	&=\sum_{k=0}^{n}{n\brace k}_{Y,\lambda}x^{k}\Gamma(k+1)=\sum_{k=0}^{n}{n\brace k}_{Y,\lambda}x^{k}k! \nonumber \\
	&=F_{n,\lambda}^{Y}(x),\quad (n\ge 0). \nonumber
\end{align}
Thus, from \eqref{27}, we obtain the following theorem.
\begin{theorem}
	For $n\ge 0$, we have
	\begin{displaymath}
		\int_{0}^{\infty}\phi_{n,\lambda}^{Y}(xy)e^{-y}dy= F_{n,\lambda}^{Y}(x).
	\end{displaymath}
\end{theorem}
From \eqref{20}, we note that
\begin{align}
\sum_{n=k}^{\infty}{n\brace k}_{Y,\lambda}\frac{t^{n}}{n!}&=\frac{1}{k!}\Big(E[e_{\lambda}^{Y}(t)]-1\Big)^{k}=\frac{1}{k!}\bigg(\sum_{i=1}^{\infty}E[(Y)_{i,\lambda}]\frac{t^{i}}{i!}\bigg)^{k} \label{28} \\
&=\sum_{n=k}^{\infty}B_{n,k}\Big(E[(Y)_{1,\lambda}],E[(Y)_{2,\lambda}],\cdots,E[(Y)_{n-k+1,\lambda}]\Big)\frac{t^n}{n!}. \nonumber	
\end{align}
Thus, by \eqref{28}, we get
\begin{equation}
{n \brace k}_{Y,\lambda}=B_{n,k}\Big(E[(Y)_{1,\lambda}],E[(Y)_{2,\lambda}],\cdots,E[(Y)_{n-k+1,\lambda}]\Big),\quad (n\ge k\ge 0). \label{29}
\end{equation}
Hence
\begin{equation}
\phi_{n,\lambda}^{Y}(y)=\sum_{k=0}^{n}{n\brace k}_{Y,\lambda}y^{k}=\sum_{k=0}^{n}	B_{n,k}\Big(E[(Y)_{1,\lambda}],E[(Y)_{2,\lambda}],\cdots,E[(Y)_{n-k+1,\lambda}]\Big)y^{k}. \label{30}
\end{equation}
By \eqref{27} and \eqref{30}, we get
\begin{align}
F_{n,\lambda}^{Y}&=\int_{0}^{\infty}\sum_{k=0}^{n} B_{n,k}\Big(E[(Y)_{1,\lambda}],E[(Y)_{2,\lambda}],\cdots,E[(Y)_{n-k+1,\lambda}]\Big)\int_{0}^{\infty}y^{k}e^{-y}dy \label{31} \\
&=\sum_{k=0}^{n}k! B_{n,k}\Big(E[(Y)_{1,\lambda}],E[(Y)_{2,\lambda}],\cdots,E[(Y)_{n-k+1,\lambda}]\Big),\quad (n\ge 0).\nonumber
\end{align}
Therefore, by \eqref{31}, we obtain the following theorem.
\begin{theorem}
	For $n\ge 0$, we have
	\begin{displaymath}
		F_{n,\lambda}^{Y}=\sum_{k=0}^{n}k! B_{n,k}\Big(E[(Y)_{1,\lambda}],E[(Y)_{2,\lambda}],\cdots,E[(Y)_{n-k+1,\lambda}]\Big).
	\end{displaymath}
\end{theorem}
For $r\in\mathbb{N}$, the {\it{probabilistic degenerate Fubini polynomials of order $r$ associated with random variable $Y$}} are defined by
\begin{equation}
\bigg(\frac{1}{1-x(E[e_{\lambda}^{Y}(t)]-1)}\bigg)^{r}=\sum_{n=0}^{\infty}F_{n,\lambda}^{(r,Y)}(x)\frac{t^{n}}{n!}.\label{32}	
\end{equation}
When $Y=1$, we have $F_{n,\lambda}^{(r,Y)}(x)=F_{n,\lambda}^{(r)}(x),\ (n\ge 0)$. \par
From \eqref{20} and \eqref{32}, we note that
\begin{align}
&\bigg(\frac{1}{1-x(E[e_{\lambda}^{Y}(t)]-1)}\bigg)^{r}=\sum_{i=0}^{\infty}\binom{-r}{i}(-1)^{i}x^{i}\Big(E[e_{\lambda}^{Y}(t)]-1\Big)^{i} \label{33} \\
&=\sum_{i=0}^{\infty}\binom{r+i-1}{i}i!x^{i}\frac{1}{i!}\Big(E[e_{\lambda}^{Y}(t)]-1\Big)^{i}=\sum_{i=0}^{\infty}\binom{r+i-1}{i}i!x^{i}\sum_{n=i}^{\infty}{n\brace i}_{Y,\lambda}\frac{t^{n}}{n!}\nonumber \\
&=\sum_{n=0}^{\infty}\sum_{i=0}^{n}\binom{r+i-1}{i}i!x^{i}{n\brace i}_{Y,\lambda}\frac{t^{n}}{n!}.\nonumber
\end{align}
Therefore, by \eqref{32} and \eqref{33}, we obtain the following theorem.
\begin{theorem}
For $n\ge 0$, we have
\begin{equation*}
F_{n,\lambda}^{(r,Y)}(x)=\sum_{i=0}^{n}\binom{r+i-1}{i}i!{n\brace i}_{Y,\lambda}x^{i}. 	
\end{equation*}	
\end{theorem}
By \eqref{23}, we get
\begin{align}
&\sum_{n=1}^{\infty}F_{n,\lambda}^{Y}(x)\frac{t^{n}}{n!}=\frac{1}{1-x(E[e_{\lambda}^{Y}(t)]-1)}-1=\frac{x(E[e_{\lambda}^{Y}(t)]-1)}{1-x(E[e_{\lambda}^{Y}(t)]-1)}	\label{35} \\
&=\frac{xE[e_{\lambda}^{Y}(t)]}{1-x(E[e_{\lambda}^{Y}(t)]-1)}-\frac{x}{1-x(E[e_{\lambda}^{Y}(t)]-1)}\nonumber \\
&=x\sum_{k=0}^{\infty}E[(Y)_{k,\lambda}]\frac{t^{k}}{k!}\sum_{l=0}^{\infty}F_{l,\lambda}^{Y}(x)\frac{t^{l}}{l!}-x\sum_{n=0}^{\infty}F_{n,\lambda}^{Y}(x)\frac{t^{n}}{n!}. \nonumber\\
&=x\sum_{k=1}^{\infty}E[(Y)_{k,\lambda}]\frac{t^{k}}{k!}\sum_{l=0}^{\infty}F_{l,\lambda}^{Y}(x)\frac{t^{l}}{l!}\nonumber \\
&=\sum_{n=1}^{\infty}x\sum_{k=1}^{n}\binom{n}{k}E\big[(Y)_{k,\lambda}\big]F_{n-k,\lambda}^{Y}(x)\frac{t^{n}}{n!}.\nonumber
\end{align}
Therefore, by comparing the coefficients on both sides of \eqref{35}, we obtain the following theorem.
\begin{theorem}
For $n\ge 1$, we have
\begin{displaymath}
F_{n,}^{Y}(x)=x\sum_{k=1}^{n}\binom{n}{k}E[(Y)_{k,\lambda}]F_{n-k,\lambda}^{Y}(x).
\end{displaymath}
\end{theorem}
From \eqref{23}, we note that
\begin{align}
&\sum_{n=0}^{\infty}F_{n+1,\lambda}^{Y}(x)\frac{t^{n}}{n!}=\frac{d}{dt}\sum_{n=0}^{\infty}F_{n,,\lambda}^{Y}(x)\frac{t^{n}}{n!}=\frac{d}{dt}\bigg(\frac{1}{1-x(E[e_{\lambda}^{Y}(t)]-1)}\bigg)\label{36} \\
&=\frac{xE[Ye_{\lambda}^{Y-\lambda}(t)]}{(1-x(E[e_{\lambda}^{Y}(t)]-1))^{2}}=\frac{x}{1-x(E[e_{\lambda}^{Y}(t)]-1)}\frac{E[Ye_{\lambda}^{Y-\lambda}(t)]}{1-x(E[e_{\lambda}^{Y}(t)]-1)}	\nonumber \\
&=x\sum_{i=0}^{\infty}F_{i,\lambda}^{Y}(x)\frac{t^{i}}{i!}\sum_{j=0}^{\infty}F_{j,\lambda}^{Y}(x)\frac{t^{j}}{j!}\sum_{m=0}^{\infty}E[Y(Y-\lambda)_{m,\lambda}]\frac{t^{m}}{m!}\nonumber \\
&=x\sum_{k=0}^{\infty}\bigg(\sum_{i=0}^{k}\binom{k}{i}F_{i,\lambda}^{Y}(x)F_{k-i,\lambda}^{Y}(x)\bigg)\frac{t^{k}}{k!}\sum_{m=0}^{\infty}E[(Y)_{m+1,\lambda}]\frac{t^{m}}{m!}\nonumber \\
&=\sum_{n=0}^{\infty}x\sum_{k=0}^{n}\sum_{i=0}^{k}\binom{k}{i}\binom{n}{k}F_{i,\lambda}^{Y}(x)F_{k-i,\lambda}^{Y}(x)E[(Y)_{n-k+1,\lambda}]\frac{t^{n}}{n!}.\nonumber
\end{align}
Therefore, by comparing the coefficients on both sides of \eqref{26}, we obtain the following theorem.
\begin{theorem}
For $n\ge 0$, we have
\begin{displaymath}
F_{n+1,\lambda}^{Y}(x)= x\sum_{k=0}^{n}\sum_{i=0}^{k}\binom{k}{i}\binom{n}{k}F_{i,\lambda}^{Y}(x)F_{k-i,\lambda}^{Y}(x)E[(Y)_{n-k+1,\lambda}].
\end{displaymath}
\end{theorem}
Now, we observe that

\begin{align}
	&\sum_{n=0}^{\infty}\frac{d^{r}}{dx^{r}}F_{n,\lambda}^{Y}(x)\frac{t^{n}}{n!}=\frac{d^{r}}{dx^{r}}\bigg(\frac{1}{1-x(E[e_{\lambda}^{Y}(t)]-1)}\bigg)=r!\frac{\big(E[e_{\lambda}^{Y}(t)]\big)^{r}}{\big(1-x(E[e_{\lambda}^{Y}(t)]-1)\big)^{r+1}}\label{37} \\
	&=r!\sum_{i=0}^{\infty}F_{i,\lambda}^{(r+1,Y)}(x)\frac{t^{i}}{i!}\sum_{k=0}^{\infty}E\big[(Y_{1}+Y_{2}+\cdots+Y_{r})_{k,\lambda}\big]\frac{t^{k}}{k!}\nonumber \\
	&=\sum_{n=0}^{\infty}\bigg(r!\sum_{i=0}^{n}F_{i,\lambda}^{(r+1,Y)}(x)E[(S_{r})_{n-i,\lambda}]\binom{n}{i}\bigg)\frac{t^{n}}{n!}.\nonumber
\end{align}
Therefore, by \eqref{37}, we obtain the following theorem.
\begin{theorem}
For $r,n\ge 0$, we have
\begin{displaymath}
\frac{d^{r}}{dx^{r}}F_{n,\lambda}^{Y}(x)= r!\sum_{i=0}^{n}F_{i,\lambda}^{(r+1,Y)}(x)E[(S_{r})_{n-i,\lambda}]\binom{n}{i}.
\end{displaymath}
\end{theorem}
From \eqref{19}, we note that
\begin{align}
&\sum_{n=0}^{\infty}\bigg(\sum_{i=0}^{\infty}E[(S_{i})_{n,\lambda}]x^{i}\bigg)\frac{t^{n}}{n!}=\sum_{i=0}^{\infty}x^{i}\sum_{n=0}^{\infty}E[(S_{i})_{n,\lambda}]\frac{t^{n}}{n!} \label{38} \\
&=\sum_{i=0}^{\infty}x^{i}E[e_{\lambda}^{S_{i}}(t)]=\sum_{i=0}^{\infty}x^{i}\Big(E[e_{\lambda}^{Y}(t)]\Big)^{i} \nonumber \\
&=\frac{1}{1-xE[e_{\lambda}^{Y}(t)]}=\frac{1}{1-x}\frac{1}{1-\frac{x}{1-x}\big(E[e_{\lambda}^{Y}(t)]-1\big)}\nonumber\\
&=\frac{1}{1-x}\sum_{n=0}^{\infty}F_{n,\lambda}^{Y}\bigg(\frac{x}{1-x}\bigg)\frac{t^{n}}{n!}. \nonumber
\end{align}
Therefore, by comparing the coefficients on both sides of \eqref{38}, we obtain the following theorem.
\begin{theorem}
For $n\ge 0$, we have
\begin{displaymath}
\frac{1}{1-x}F_{n,\lambda}^{Y}\bigg(\frac{x}{1-x}\bigg)=\sum_{i=0}^{\infty}E[(S_{i})_{n,\lambda}]x^{i}.
\end{displaymath}
\end{theorem}
Taking $x=\frac{1}{2}$, we get
\begin{displaymath}
\sum_{i=0}^{\infty}E[(S_{i})_{n,\lambda}]\bigg(\frac{1}{2}\bigg)^{i}=2F_{n,\lambda}^{Y},\quad (n\ge 0).
\end{displaymath}
Let $Y$ be the Poisson random variable with parameter $\alpha(>0)$. Then we have
\begin{equation}
E[e_{\lambda}^{Y}(t)]=\sum_{n=0}^{\infty}e_{\lambda}^{n}(t)\frac{\alpha^{n}}{n!}e^{-\alpha}=e^{\alpha(e_{\lambda}(t)-1)}. \label{39}	
\end{equation}
From \eqref{39}, we have
\begin{align}
\sum_{n=0}^{\infty}F_{n,\lambda}^{Y}(x)\frac{t^{n}}{n!}&=\frac{1}{1-x(E[e_{\lambda}^{Y}(t)]-1)}=\frac{1}{1-x(e^{\alpha(e_{\lambda}(t)-1)}-1)}\label{40} \\
&=\sum_{i=0}^{\infty}F_{i}(x)\alpha^{i}\frac{1}{i!}\big(e_{\lambda}(t)-1\big)^{i}=\sum_{i=0}^{\infty}F_{i}(x)\alpha^{i}\sum_{n=i}^{\infty}{n \brace i}_{\lambda}\frac{t^{n}}{n!} \nonumber \\
&=\sum_{n=0}^{\infty}\bigg(\sum_{i=0}^{n}F_{i}(x){n\brace i}_{\lambda}\alpha^{i}\bigg)\frac{t^{n}}{n!},\nonumber	
\end{align}
where $F_{i}(x)$ are the Fubini polynomials given by $\frac{1}{1-x(e^t -1)}=\sum_{i=0}^{\infty}F_{i}(x)\frac{t^i}{i!}$.
Therefore, by \eqref{40}, we obtain the following theorem.
\begin{theorem}
Let $Y$ be the Poisson random variable with parameter $\alpha(>0)$. Then we have
\begin{displaymath}
F_{n,\lambda}^{Y}(x)=\sum_{i=0}^{n}F_{i}(x){n \brace i}_{\lambda}\alpha^{i},\quad (n\ge 0).
\end{displaymath}
\end{theorem}
Let $Y$ be the Poisson random variable with parameter $\alpha>0$. Then we have
\begin{equation}
\Big(E[e_{\lambda}^{Y}(t)]\Big)^{k}=e^{k\alpha(e_{\lambda}(t)-1)}=\sum_{n=0}^{\infty}\phi_{n,\lambda}(k\alpha)\frac{t^{n}}{n!}.\label{41}	
\end{equation}
and
\begin{equation}
\Big(E[e_{\lambda}^{Y}(t)]\Big)^{k}
= \sum_{n=0}^{\infty}E[(S_{k})_{n,\lambda}]\frac{t^{n}}{n!}. \label{42}
\end{equation}
Thus, by \eqref{41} and \eqref{42}, we get
\begin{equation}
E\big[(S_{k})_{n,\lambda}\big]=\phi_{n,\lambda}(k\alpha),\quad (n\ge 0). \label{43}
\end{equation}
From Theorem 2.10 and \eqref{43}, we have
\begin{equation}
\sum_{k=0}^{\infty}\phi_{n,\lambda}(k\alpha)x^{k}=\sum_{k=0}^{\infty}E\Big[(S_{k})_{n,\lambda}\Big]x^{k}=\frac{1}{1-x}F_{n,\lambda}^{Y}\bigg(\frac{x}{1-x}\bigg).\label{44}	
\end{equation}
Therefore, by \eqref{44}, we obtain the following theorem.
\begin{theorem}
Let $Y$ be the Poisson random variable with parameter $\alpha(>0)$. For $n\ge 0$, we have
\begin{displaymath}
\sum_{k=0}^{\infty}\phi_{n,\lambda}(k\alpha)x^{k}=\sum_{k=0}^{\infty}E\Big[(S_{k})_{n,\lambda}\Big]x^{k}=\frac{1}{1-x}F_{n,\lambda}^{Y}\bigg(\frac{x}{1-x}\bigg).
\end{displaymath}
\end{theorem}
By using Theorem 2.6 and \eqref{19}, we note that
\begin{align}
\bigg(\frac{1}{1-x}\bigg)^{r+1}F_{n,\lambda}^{(r+1,Y)}\bigg(\frac{x}{1-x}\bigg)&=\bigg(\frac{1}{1-x}\bigg)^{r+1}\sum_{l=0}^{n}{n \brace l}_{Y,\lambda}\binom{l+r}{l}l!\bigg(\frac{x}{1-x}\bigg)^{l}\label{45} \\
&=\sum_{l=0}^{n}{n \brace l}_{Y,\lambda}l!\binom{l+r}{l}x^{l}\bigg(\frac{1}{1-x}\bigg)^{l+r+1}\nonumber \\
&=\sum_{l=0}^{n}\sum_{k=0}^{\infty}{n \brace l}_{Y,\lambda}l!\binom{l+r}{l}\binom{k+l+r}{k}x^{k+l}\nonumber \\
&=\sum_{l=0}^{n}\bigg(\sum_{k=l}^{\infty}{n \brace l}_{Y,\lambda}l!\binom{l+r}{l}\binom{k+r}{k-l}\bigg)x^{k}\nonumber \\
&=\sum_{k=0}^{n}\binom{k+r}{k}x^{k} \sum_{l=0}^{k}{n \brace l}_{Y,\lambda}(k)_{l}+\sum_{k=n+1}^{\infty}\binom{k+r}{k}x^{k} \sum_{l=0}^{n}{n \brace l}_{Y,\lambda}(k)_{l}\nonumber \\
&=\sum_{k=0}^{n}\binom{k+r}{k}x^{k} \sum_{l=0}^{k}{n \brace l}_{Y,\lambda}(k)_{l}+\sum_{k=n+1}^{\infty}\binom{k+r}{k}x^{k} \sum_{l=0}^{k}{n \brace l}_{Y,\lambda}(k)_{l}\nonumber \\
&=\sum_{k=0}^{\infty}\binom{k+r}{k}x^{k} \sum_{l=0}^{k}{n \brace l}_{Y,\lambda}(k)_{l}=\sum_{k=0}^{\infty}\binom{k+r}{k}x^{k}E\Big[(S_{k})_{n,\lambda}\Big].	\nonumber
\end{align}
Therefore, by \eqref{45}, we obtain the following theorem.
\begin{theorem}
For $n\ge 0$, we have
\begin{displaymath}
\bigg(\frac{1}{1-x}\bigg)^{r+1}F_{n,\lambda}^{(r+1,Y)}\bigg(\frac{x}{1-x}\bigg)= \sum_{k=0}^{\infty}\binom{k+r}{k}x^{k}E\Big[(S_{k})_{n,\lambda}\Big].
\end{displaymath}
\end{theorem}
When $Y=1$, we have
\begin{displaymath}
\bigg(\frac{1}{1-x}\bigg)^{r+1}F_{n,\lambda}^{(r+1)}\bigg(\frac{x}{1-x}\bigg)=\sum_{k=0}^{\infty}\binom{k+r}{k}x^{k}(k)_{n,\lambda}.
\end{displaymath}
Now, we observe from \eqref{21} and Theorem 2.6 that
\begin{align}
&\int_{0}^{\infty}y^{r-1}\phi_{n,\lambda}^{Y}(xy)e^{-y}dy=\sum_{k=0}^{n}{n \brace k}_{Y,\lambda}x^{k}\int_{0}^{\infty}y^{r+k-1}e^{-y}dy \label{46}\\
&=\sum_{k=0}^{n}{n \brace k}_{Y,\lambda}x^{k}\Gamma(r+k)=\Gamma(r)\sum_{k=0}^{n}{n\brace k}_{Y,\lambda}x^{k}\binom{r+k-1}{k}k!\nonumber \\
&=\Gamma(r)F_{n,\lambda}^{(r,Y)}(x),\quad (r\in\mathbb{N}).\nonumber
\end{align}
Therefore, by \eqref{46}, we obtain the following theorem.
\begin{theorem}
For $n\ge 0$ and $r\ge 1$, we have
\begin{equation*}
F_{n,\lambda}^{(r,Y)}(x)=\frac{1}{\Gamma(r)}\int_{0}^{\infty}y^{r-1}\phi_{n,\lambda}^{Y}(xy)e^{-y}dy.
\end{equation*}
\end{theorem}
From \eqref{32} and Theorem 2.14, we have
\begin{align}
&\bigg(\frac{1}{1-x(E[e_{\lambda}^{Y}(t)]-1)}\bigg)^{r}=\frac{1}{\Gamma(r)}\int_{0}^{\infty}y^{r-1}\sum_{n=0}^{\infty}\phi_{n,\lambda}^{Y}(xy)\frac{t^n}{n!}e^{-y}dy\label{47} \\
&=\frac{1}{\Gamma(r)}\int_{0}^{\infty}y^{r-1}e^{xy(E[e_{\lambda}^{Y}(t)]-1)}e^{-y}dy \nonumber \\
&=\frac{1}{\Gamma(r)}\int_{0}^{\infty} y^{r-1}e^{y(xE[e_{\lambda}^{Y}(t)]-1-x)}dy \nonumber.
\end{align}
By \eqref{47}, we get
\begin{displaymath}
\frac{1}{\Gamma(r)}\int_{0}^{\infty} y^{r-1}e^{y(xE[e_{\lambda}^{Y}(t)]-1-x)}dy= \bigg(\frac{1}{1-x(E[e_{\lambda}^{Y}(t)]-1)}\bigg)^{r},
\end{displaymath}
where $r$ is a positive integer. \par
We observe that
\begin{align}
\sum_{n=0}^{\infty}F_{n+1}^{(r,Y)}(x)\frac{t^{n}}{n!}&=\frac{d}{dt}\sum_{n=0}^{\infty}F_{n,\lambda}^{Y}(x)\frac{t^{n}}{n!} \label{48} \\
&=\frac{d}{dt}\bigg(\frac{1}{1-x(E[e_{\lambda}^{Y}(t)]-1)}\bigg)^{r}=\frac{rxE[e_{\lambda}^{Y-\lambda}(t)Y]}{\big(1-x(E[e_{\lambda}^{Y}(t)]-1)\big)^{r+1}}\nonumber \\
&= \frac{rx}{\big(1-x(E[e_{\lambda}^{Y}(t)]-1)\big)^{r}} \frac{1}{1-x(E[e_{\lambda}^{Y}(t)]-1)}E[e_{\lambda}^{Y-\lambda}(t)Y]\nonumber \\
&=rx\sum_{l=0}^{\infty}F_{l,\lambda}^{(r,Y)}(x)\frac{t^{l}}{l!}\sum_{i=0}^{\infty}F_{i,\lambda}^{Y}(x)\frac{t^{i}}{i!}\sum_{j=0}^{\infty}E[(Y)_{j+1,\lambda}]\frac{t^{j}}{j!} \nonumber \\
&=rx\sum_{l=0}^{\infty}F_{l,\lambda}^{(r,Y)}(x)\frac{t^{l}}{l!}\sum_{k=0}^{\infty}\bigg(\sum_{j=0}^{k}\binom{k}{j}F_{k-j,\lambda}^{Y}(x)E[(Y)_{j+1,\lambda}]\bigg)\frac{t^{k}}{k!} \nonumber\\
&=\sum_{n=0}^{\infty}rx\sum_{k=0}^{n}\sum_{j=0}^{k}\binom{n}{k}\binom{k}{j}F_{n-k,\lambda}^{(r,Y)}(x)F_{k-j,\lambda}^{Y}(x)E[(Y)_{j+1,\lambda}]\frac{t^{n}}{n!}. \nonumber
\end{align}
Therefore, by \eqref{48}, we obtain the following theorem.
\begin{theorem}
For $n\ge 0$, we have
\begin{displaymath}
F_{n+1,\lambda}^{(r,Y)}(x)= rx\sum_{k=0}^{n}\sum_{j=0}^{k}\binom{n}{k}\binom{k}{j}F_{n-k,\lambda}^{(r,Y)}(x)F_{k-j,\lambda}^{Y}(x)E[(Y)_{j+1,\lambda}].
\end{displaymath}
\end{theorem}
Let $Y$ be the Bernoulli random variable with probability of success $p$. Then we have
\begin{equation}
E[e_{\lambda}^{Y}(t)]=\sum_{k=0}^{1}e_{\lambda}^{k}(t)p(k)=1-p+pe_{\lambda}(t)=1+p(e_{\lambda}(t)-1). \label{49}	
\end{equation}
By \eqref{49}, we get
\begin{align}
\sum_{n=0}^{\infty}F_{n,\lambda}^{Y}(x)\frac{t^{n}}{n!}&=\frac{1}{1-x(E[e_{\lambda}^{Y}(t)]-1)}=\frac{1}{1-xp(e_{\lambda}(t)-1)} \label{50}\\
&=\sum_{n=0}^{\infty}F_{n,\lambda}(xp)\frac{t^{n}}{n!}.\nonumber
\end{align}
Therefore, by comparing the coefficients on both sides of \eqref{50}, we obtain the following theorem.
\begin{theorem}
Let $Y$ be the Bernoulli random variable with probability of success $p$. For $n\ge 0$, we have
\begin{displaymath}
F_{n,\lambda}^{Y}(x)=F_{n,\lambda}(xp).
\end{displaymath}
\end{theorem}
\section{Conclusion}
In this paper, we studied by using generating functions the probabilistic degenerate Fubini polynomials associated with $Y$ and the probabilistic degenerate Fubini polynomials of order $r$ associated with $Y$, as probabilistic versions of the degenerate Fubini polynomials and the degenerate Fubini polynomials of order $r$, respectively. Here $Y$ is a random variable such that the moment generating function of $Y$ exists in a neighborhood of the origin. In more detail, we derived several explicit  expressions of $F_{n, \lambda}^{Y}(x)$ (see Theorems 2.1, 2.2, 2.4) and those of $F_{n, \lambda}^{r,Y}(x)$ (see Theorems 2.6, 2.14). We obtained a recurrence relations for $F_{n, \lambda}^{Y}(x)$ (see Theorem 2.7), and another one (see Theorem 2.8) together with its generalization (see Theorem 2.15). We expressed the $r$th derivative of $F_{n, \lambda}^{Y}(x)$ in terms of $F_{i,\lambda}^{(r+1,Y)}(x)$ (see Theorem 2.9). We showed the identity $\frac{1}{1-x}F_{n,\lambda}^{Y}\bigg(\frac{x}{1-x}\bigg)=\sum_{i=0}^{\infty}E[(S_{i})_{n,\lambda}]x^{i}$ (see Theorem 2.10) and its generalization (see Theorem 2.13). We deduced an explicit expression for $F_{n, \lambda}^{Y}(x)$ when $Y \sim \Gamma(1,1)$ (see Theorem 2.3) and also that when $Y$ is the Poisson random variable with parameter $\alpha$ (see Theorem 2.11). We proved that $\frac{1}{1-x}F_{n,\lambda}^{Y}\bigg(\frac{x}{1-x}\bigg)=\sum_{k=0}^{\infty}\phi_{n,\lambda}(k\alpha)x^{k}$ when $Y$ is the Poisson random variable with parameter $\alpha$ (see Theorem 2.12). We showed $F_{n,\lambda}^{Y}(x)=F_{n,\lambda}(xp)$ when
$Y$ be the Bernoulli random variable with probability of success $p$ (see Theorem 2.16). \par
As one of our future projects, we would like to continue to study degenerate versions, $\lambda$-analogues and probabilistic versions of many special polynomials and numbers and to find their applications to physics, science and engineering as well as to mathematics.

\section*{Declarations}

{\bf{Ethical Statement}} The submitted article is an original research paper and has not been published anywhere
else. There are no clinical trials, animal research, or human trials included in this paper.

{\bf{Funding}}
Tis research was funded by the National Natural Science Foundation of China (No.12271320) and Key Research and Development Program of Shaanxi (No. 2023
ZDLGY-02).

{\bf { Data Availability}} There is no data used in this article.

{\bf{Conflict of Interest}} The authors declare no competing interests.

\end{document}